\documentclass{conm-p-l}
\usepackage{amssymb,graphicx}
\usepackage{url}
\usepackage[OT2,OT1]{fontenc}
\DeclareFontEncoding{OT2}{}{}
\newcommand\cy[1]{{\fontencoding{OT2}\selectfont #1}}
\newtheorem{theorem}{Theorem}[section]
\newtheorem{lemma}[theorem]{Lemma}
\newtheorem{corollary}[theorem]{Corollary}

\newtheorem{definition}[theorem]{Definition}

\newtheorem{conjecture}[theorem]{Conjecture}

\newtheorem{problem}[theorem]{Problem}
\newtheorem{ComplexityConjecture}[theorem]{Hyperbolic Complexity
Conjecture} \theoremstyle{definition}
 \theoremstyle{definition}
\newtheorem{Mom4Conjecture}[theorem]{Mom-$4$ Conjecture} \theoremstyle{definition}
\newtheorem{remark}[theorem]{Remark}
\newtheorem{remarks}[theorem]{Remarks}
\newcommand{\Vol}{\operatorname{Vol}}
\newcommand{\Area}{\operatorname{Area}}
\newcommand{\BH}{\mathbb{H}}
\newcommand{\BR}{\mathbb{R}}
\newcommand{\BC}{\mathbb{C}}

\newcommand{\BZ}{\mathbb{Z}}
\newcommand{\BQ}{\mathbb{Q}}
\newcommand{\mr}[1]{\mathrm{#1}}

\newcommand{\mc}[1]{\mathcal{#1}}
\newcommand{\barM}{\bar{M}}

\begin{document}
\title{Mom technology and hyperbolic 3-manifolds}
\author{David Gabai}\thanks{The first author was partially supported by NSF grants
DMS-0554374 and DMS-0504110.}\address{Department of
Mathematics\\Princeton University\\Princeton, NJ}
\author{Robert Meyerhoff}\thanks{The second author was partially supported by NSF grants
DMS-0553787 and DMS-0204311.}\address{Department of Mathematics\\Boston
College\\Chestnut Hill, MA}
\author{Peter Milley}\thanks{The third author was partially supported by NSF grant
DMS-0554624 and ARC Discovery grant DP0663399.}\address{Department
of Mathematics and Statistics\\University of Melbourne\\Melbourne,
Australia}


\date{}
\maketitle

\section{Introduction}

This is an expository paper on the work of the authors, found in
\cite{GMM2}, \cite{GMM3}, and \cite{M}, which proves that the Weeks
manifold is the unique closed orientable $3$-manifold of minimum volume
and determines the ten one-cusped hyperbolic $3$-manifolds of volume at
most $2.848$.  Our work has focused on \emph{Mom technology}, which
has proven to be effective in determining these low-volume manifolds
and has the potential for vast generalization.

This introduction will survey a selection of results on volumes of
hyperbolic 3-manifolds.  The body of the paper will outline the recent
work of the authors mentioned above, and the last section will discuss
outstanding open problems in this area.

Unless stated otherwise, all manifolds in this paper will be
orientable.

The theory of volumes of hyperbolic $3$-manifolds has received
tremendous interest over the last 40 years.  In 1968 Mostow proved via
his famous Rigidity theorem in \cite{Most} that volume is a topological
invariant of closed manifolds. This was extended to finite-volume
hyperbolic $3$-manifolds by Marden \cite{Ma} and Prasad \cite{Pra}.
In the mid 1970's, Troels Jorgensen proved
that for any constant $C$ the collection of all complete hyperbolic
$3$-manifolds of volume at most $C$ is obtained from a finite
collection of cusped manifolds using the operation of Dehn filling;
see \cite{Th1} or \cite{Gr}.
Using these ideas together with those of Gromov, Thurston showed in
1977 that the subset $V\subset\BR$ of volumes of complete hyperbolic
$3$-manifolds is a closed well-ordered set of order type
$\omega^\omega$.  Furthermore, there are only finitely many manifolds
of a given volume.  Thurston also showed that every $n$-cusped
hyperbolic $3$-manifold is a limit volume of hyperbolic $3$-manifolds
with $n-1$ cusps and that any filling on an $n$-cusped manifold yields
a manifold of lower volume.  In particular, the smallest volume
3-manifold is closed.

We remind the reader that saying $V$ is of order type $\omega^\omega$
implies that there is a smallest element $v_1$, a next smallest element
$v_2$, and so on, with $v_\omega$ being the first limit element, then
$v_{\omega+1}$, $v_{\omega+2}$, \ldots\ limiting on $v_{2\omega}$, the
second limit volume.  The limit volumes $v_{\omega}$, $v_{2\omega}$,
$v_{3\omega}$, \ldots\ limit on $v_{\omega^2}$, the first limit of
limit volumes.  The volumes $v_{n\omega}$ represent the volumes of $1$-cusped
manifolds, $v_{m\omega^2}$ the volumes of $2$-cusped manifolds, and so
on (although there may be compact manifolds or manifolds with fewer cusps with these volumes as
well). In contrast, for $n\ge 4$ the set of volumes of complete
hyperbolic $n$-manifolds is discrete; see \cite{Wa}.

\begin{theorem}
The Matveev-Fomenko-Weeks manifold is the unique smallest volume
closed hyperbolic 3-manifold.  In particular $v_1=0.9427\ldots$.
\label{main}
\end{theorem}

This manifold, commonly referred to as the Weeks manifold for short,
can be constructed as the $(5,1)$, $(5,2)$ Dehn filling on the
complement of the Whitehead link. It was discovered and its volume was
computed independently by Matveev and Fomenko in \cite{MF} and
Weeks using his SnapPea program \cite{W}.  Matveev and Fomenko were
motivated by Novikov who was interested in volumes for physical and
dynamical reasons; see \cite{NS}.  Independently, Przytycki
asked whether this manifold was the smallest closed manifold.  See
\cite{BPZ} for related questions.

There were many partial results towards this problem: \cite{Mey1},
\cite{Mey2}, \cite{GM1}, \cite{GMT} (using \cite{GM2}), \cite{Prz1},
\cite{MM1}, \cite{A1},\cite{Prz2}, and (Agol - Dunfield)
\cite{AST}. The ideas of many of these results were used in our proof
of Theorem \ref{main}. Their role will be discussed further in the
body of the paper.  In addition, Weeks' remarkable SnapPea
program was indispensable for providing experimental data and
geometric intuition.

Other results along these lines include: the smallest cusped
hyperbolic orbifold in \cite{Mey3}, the smallest cusped manifold
(orientable or not) in \cite{Ada}, the smallest cusped orientable
manifold in \cite{CM}, the smallest arithmetic hyperbolic $3$-orbifold
in \cite{CF}, the smallest compact manifold with totally geodesic
boundary in \cite{KM}, the smallest $2$-cusped hyperbolic
$3$-manifolds in \cite{A2}, and via a tour de force the smallest
hyperbolic orbifold in \cite{GM3} and \cite{MM3}.

Underlying some of these results are the very useful packing results
of \cite{Bor}, \cite{Prz2}, \cite{Prz3}, and \cite{Miy}.  In particular Miyamoto \cite{Miy} gives remarkable results on manifolds with totally geodesic boundary.

In other directions there is the long series of papers by Culler,
Shalen and their co-authors which give lower bounds on volumes for
hyperbolic $3$-manifolds that satisfy certain topological constraints.
See for example \cite{CHS1}, \cite{CHS2}, \cite{ACS1}, \cite{ACS2}, and
their references.

See Milnor's paper \cite{Mil} for a detailed history of hyperbolic
geometry through about 1980.  Its appendix contains volume formulae
for ideal tetrahedra and in particular a proof that the regular ideal
simplex is the one of maximal volume.  This last result was proven in
higher dimensions by Haggerup and Munkholm \cite{HM}.

\smallskip
\noindent\emph{Acknowledgments:} We thank Ian Agol and Nathan Dunfield for
their insightful
comments on our problem list and Nathan Dunfield for also sharing many of his
experimental results for this section.  We are also grateful to Alan Reid 
and Walter Neumann for contributing several number-theoretic problems. 

\section{Outline Of The Proof}

For basic facts about hyperbolic 3-manifolds see \cite{Th1},
\cite{Ra}, or \cite{BP}.

By a \emph{cusped hyperbolic manifold} we mean a complete non-compact
hyperbolic $3$-manifold $M$ with finite volume. By a \emph{compact
hyperbolic manifold} we mean a compact manifold whose interior supports a
complete hyperbolic metric.  Such a manifold is either closed or its
boundary is a union of tori.  Any cusped hyperbolic $3$-manifold $M$
naturally compactifies to a compact hyperbolic $3$-manifold $\barM$.  A
\emph{Dehn filling} on the cusped manifold $M$ is the interior of a
manifold obtained by attaching solid tori to various components of
$\partial\barM$.

The basic idea behind our efforts to identify low-volume
hyperbolic $3$-manifolds, particularly the smallest closed manifold,
is as follows.
Given $V\in \BR$ we identify a finite, reasonable set of cusped
manifolds $M_1$, \ldots, $M_k$ such that every hyperbolic $3$-manifold of
volume at most $V$ is obtained by filling at least one of the
$M_i$'s.  Then we identify all of the manifolds obtained by filling
the $M_i$'s that have volume at most $V$.  When $V=2.848$ and the
manifolds in question have exactly one cusp, the first step is carried
out in \cite{GMM2} and \cite{GMM3} and the last step is carried out in
\cite{M}.  We showed that exactly $10$ one-cusped manifolds, the
first $10$ in the one-cusped SnapPea census, have volume at most
$2.848$. By Agol \cite{ACS1} the smallest volume hyperbolic $3$-manifold is
obtained by filling a $1$-cusped manifold of volume at most $2.848$.
This result makes crucial use of Agol-Dunfield \cite{AST} which in turn
makes crucial use of Perelman \cite{P1,P2} and \cite{GMT}. 
Further analysis in
\cite{M} of these $10$ manifolds identifies the Weeks manifold as the
unique one with volume at most .$9428$.

Sections 3 through 6 discuss the work of \cite{GMM2} and
\cite{GMM3}. Section 7 discusses \cite{M}. Section 8 describes
previous work in identifying the smallest volume manifold as well as
methods developed in those works that play an important role in the
final resolution.  Section 9 gives more detail on some final issues,
in particular the completion of the proof of Theorem \ref{main}.  In
section 10, problems and directions for future research are presented.

For now we describe, in the context of one-cusped manifolds, the idea
behind finding the $M_i$'s cited above using the notion
of \emph{Mom technology} introduced in \cite{GMM2}.  A Mom-$n$
manifold is a compact manifold with a particular type of handle
structure, where $n$ denotes the complexity of the structure. More
formally, we have the following:


\begin{definition} A \emph{Mom-$n$} is a triple $(M,T,\Delta)$ where
\begin{itemize}
\item $M$ is a compact 3-manifold with boundary a non-empty disjoint
union of tori,
\item $T$ is a component of $\partial M$, a small neighborhood of
  which can be identified with $T\times I$, and
\item $\Delta$ is a handle structure on $M$ of the following type.
Starting with $T\times I$, where $T$ is identified with $T\times\{0\}$,
attach $n$ $1$-handles to $T\times\{1\}$ followed by $n$
$2$-handles on the ``$T\times 1$-side'' according to the following rule:
counting with multiplicity, each $1$-handle
meets at least two $2$-handles and each $2$-handle is attached to
\emph{exactly} three $1$-handles, again counting with multiplicity.
\end{itemize}
\label{defn:Mom-n}
\end{definition}

We say that a Mom-$n$ is \emph{hyperbolic} if its interior is
hyperbolic. Note that a Mom-$n$ manifold has at least two boundary
components. There are exactly $3$ hyperbolic Mom-$2$ manifolds, $18$
hyperbolic Mom-$3$ manifolds and conjecturally $117$ hyperbolic
Mom-$4$ manifolds; see \cite{GMM2}.  It is proven in \cite {GMM3} that
every $1$-cusped
hyperbolic $3$-manifold of volume at most $2.848$ is obtained by
filling a hyperbolic Mom-$n$ manifold, where $n\le 3$, using the
notion of an \emph{internal Mom-$n$ structure}:

\begin{definition}
An \emph{internal Mom-$n$ structure} on a (closed or cusped) hyperbolic
$3$-manifold $N$ consists of a Mom-$n$ $(M,T,\Delta)$ together
with an embedding $i:M\rightarrow N$ such that the image of each
component of $\partial M$ either cuts off a cusp neighbourhood
or a solid torus to the outside of $M$.
The \emph{Mom number} of $N$ is the minimal $n$ such that $N$ has an
internal Mom-$n$ structure with $M$ hyperbolic.
\label{defn:Mom-n-struct}
\end{definition}


Roughly speaking the Mom number of $N$ is the smallest $n$ such that there
exists an essentially embedded hyperbolic Mom-$n$ manifold inside of $N$.
So given a low-volume one-cusped 3-manifold, the goal is to find an
internal Mom-$n$ structure with $n$ small.  One intuitive reason for
expecting the existence of such structures is as follows. Let $W$ be
the maximal cusp of the one-cusped manifold $N$.  Topologically,
$W=T^2\times[1,\infty)$ with a finite set of pairs of points of
$\partial W$ identified with each other, and geometrically each $T^2
\times t$ is a horotorus. Further each $x\times [1,\infty)$ is an
isometrically embedded geodesic.
By straightforward geometric reasoning,
$W$ (resp. $\partial W$) already has a decent amount of volume (resp.
area). Now expand $W$ out in the standard Morse-theoretic way.
If $N$ has low volume then this forces $W$ to rapidly encounter
itself, leading
to the creation of handles of index at least $1$.  If at some moment
$n$ $1$-handles, $n$ $2$-handles, and no $3$-handles have been
created, the resulting
manifold will have Euler characteristic $0$ and hence is a
candidate for being a hyperbolic Mom-$n$ manifold.  See \S3 for
explicit examples.

Our motivation for the terminology \emph{Mom manifold} is based on
Thurston's parent-child relationship between compact 3-manifolds.
If $N$ is obtained by filling $M$, then Thurston called $M$ the
\emph{parent} and $N$ the \emph{child}.  

\begin{definition} As per Shubert and Matveev,
denote the intersections of the $1$-handles and $2$-handles with
$T\times\{1\}$ by \emph{islands} and \emph{bridges} respectively, and
the closed complement of the islands and bridges in $T\times\{1\}$ by
\emph{lakes}.  A Mom-$n$ is \emph{full} if the
lakes are all simply connected. \end{definition}

Note that if $n$ is minimal, then
being full is a necessary condition for being hyperbolic, as
otherwise either $M$ has an embedded essential annulus joining
$T\times\{1\}$ to $T\times\{0\}$ or the lake is compressible in $M$.
In that case either $M$ is reducible, contradicting hyperbolicity,
or an essential compressing disk for the lake is boundary parallel
in $M$ and hence $n$ can be reduced, contradicting minimality of $n$.

Here is a very brief outline for finding a hyperbolic internal Mom-$3$
structure in a one-cusped hyperbolic $3$-manifold $N$ of volume at
most $2.848$.  The preimage $\{B_i\}$ in $\BH^3$ of the maximal cusp
is a union of horoballs.  We consider $ \pi_1(N)$-orbits of unordered
pairs of such balls and $\pi_1(N)$-orbits of unordered triples of
horoballs.  Note that a triple of balls involves three pairs, however
two such pairs may be in equivalent classes.  A set of $n$ classes of
triples which involve exactly $n$ classes of pairs is called a
\emph{combinatorial Mom-$n$ structure}.  All one-cusped manifolds of
volume $2.848$ or less have a combinatorial Mom-$n$ structure with $n\le 3$.
See \S 4 for more details.  A combinatorial Mom-$n$ structure gives
rise to an immersed \emph{geometric} Mom-$n$ structure, i.e. the cores
of the $1$ and $2$ handles are totally geodesic.  Here, pairs of balls
give rise to $1$-handles and triples of balls give rise to
$2$-handles.  Unfortunately, the handles may intersect each other (or
even the cusp) in undesirable ways. But under the $2.848$ constraint a
``controlled'' combinatorial Mom-$n$ structure, $n\le 3$, can be
found, and after simplification this structure can be ``promoted'' to
produce an embedded geometric Mom-$k$ structure, $k\le n$.  Then if the
Mom-$k$ manifold is not hyperbolic, geometric and topological
arguments produce a hyperbolic internal Mom-$p$ structure in $N$
with $p\le k$. See \S 5 for more details.  The issue of enumerating
the Mom-$n$ manifolds, $n\le 3$ is discussed in \S 6.

\section{Example: the figure-$8$ knot and the manifold $m069$.}

As a first example, let $N$ be the complement of the figure-$8$ knot in
$S^3$; we will construct an internal Mom-$2$ structure. For the torus
$T$, we choose the boundary of a cusp neighborhood, i.e. the knot
torus. Then to $T\times I$ we add two $1$-handles and two $2$-handles
as shown in Figure
\begin{figure}
\begin{center}
\includegraphics{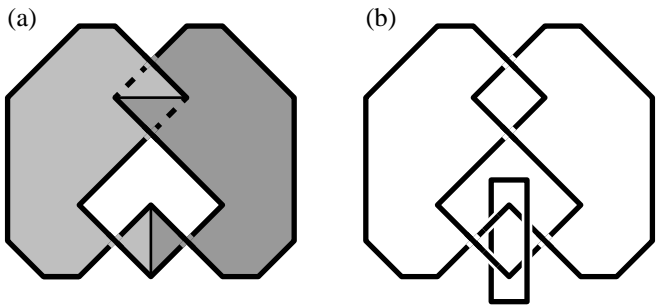}
\end{center}
\caption{(a) An internal Mom-$2$ structure on the figure-$8$ knot
  complement. (b) The interior of the corresponding submanifold $M$ is
  the complement of this Whitehead link.}
\label{fig:m004_mom}
\end{figure}
\ref{fig:m004_mom}(a), to obtain a handle structure
$\Delta$ on a submanifold $M$. Note that each $2$-handle is connected
to three $1$-handles, counting multiplicity, and the boundary of $M$
consists of two tori, one of which is $T\times\{0\}$ and the other of
which is on the ``outside'' of $M$ as seen in the figure. Thus
$(M,T,\Delta)$ is a Mom-$2$. Furthermore, the component of
$\partial M$ which is not $T\times\{0\}$ clearly bounds a solid
torus in $N\setminus M$, and $M$ is hyperbolic as discussed below.
Hence the embedding defines an internal Mom-$2$ structure on $N$.

Readers familiar with the figure-$8$ knot complement will note
immediately that the $1$-handles and $2$-handles in this example are
neighborhoods of geodesic arcs and totally geodesic surfaces. In
fact the $1$-handles and $2$-handles are neighborhoods of edges and
faces in the canonical ideal triangulation of this manifold; this
property is shared by all hyperbolic internal Mom-$n$ structures in
one-cusped
manifolds as far as the authors are aware. This makes Mom-$n$
structures easy to find in such manifolds once the
canonical ideal triangulation of the manifold is known.

The interior of the submanifold $M$ in this example is clearly homeomorphic to
$N\setminus\gamma$ where $\gamma$ is the curve indicated in Figure
\ref{fig:m004_mom}(b). The complement of this link is
homeomorphic to the complement of the Whitehead link after a Dehn
twist along $\gamma$; hence $(M,T,\Delta)$ is a hyperbolic
Mom-$2$. Another choice of internal Mom-$2$ structure on the figure-$8$
knot complement is shown in Figure
\begin{figure}
\begin{center}
\includegraphics{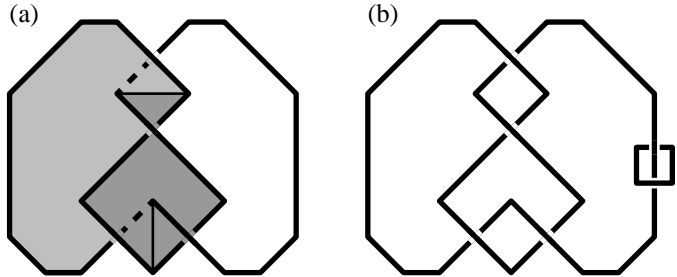}
\end{center}
\caption{(a) A non-hyperbolic internal Mom-$2$ structure on the figure-$8$ knot
  complement. (b) The corresponding non-hyperbolic link.}
\label{fig:m004_badmom}
\end{figure}
\ref{fig:m004_badmom}(a). In this case we
get a Mom-$2$ $(M',T,\Delta')$ which is not hyperbolic because it is
not full; the rightmost strand in the knot diagram forms an annular lake,
and consequently $M'$ contains an essential annulus.

The previous example, while very simple, doesn't illuminate the
Morse-theoretic nature of internal Mom-$n$ structures, nor does it illustrate
the simplest way of finding internal Mom-$n$ structures in more general
hyperbolic $3$-manifolds. So consider the one-cusped manifold known as
$m069$ in the SnapPea census. A \emph{cusp diagram} for this manifold as
produced by SnapPea is shown in Figure
\begin{figure}
\begin{center}
\includegraphics{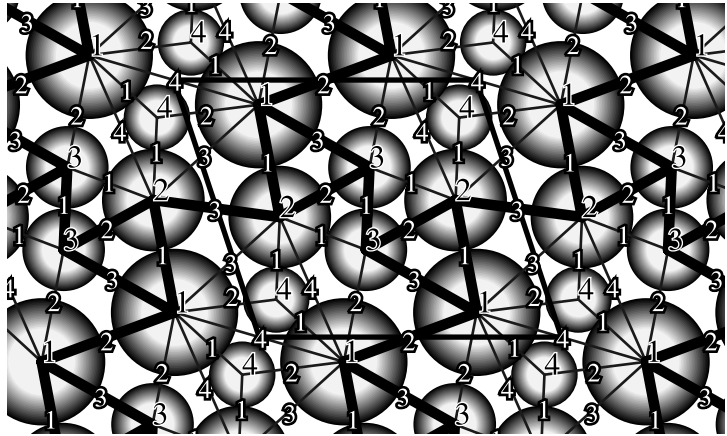}
\end{center}
\caption{A cusp diagram (from data provided by SnapPea) for the
  manifold $m069$, with an internal Mom-$3$ structure highlighted.}
\label{fig:m069}
\end{figure}
\ref{fig:m069}, including the link of the canonical ideal
triangulation of $m069$ and the shadows of nearby horoballs in the
universal cover.
Not seen is the horoball $B_\infty$ which lies above a horizontal
plane. The label on an edge in the diagram indicates the
\emph{orthopair class} of the
pair of horoballs below the endpoints of the edge, while the label
on a horoball indicates the orthopair class of the pair of horoballs
consisting of that horoball and $B_\infty$. The notion of an orthopair
class is defined in the next section; for now we merely note that
the lower the label, the shorter (possibly equal) the distance between
two horoballs of that pair.

There exists an internal Mom-$3$ structure on $m069$ where the three
$1$-handles correspond to the three shortest return paths from the
cusp torus to itself. In other words, the $1$-handles are the first
three $1$-handles created by expanding the cusp neighborhood $W$ as
described in the previous section. These three $1$-handles meet the
cusp torus in six islands, which will occur at the centers of the
horoball shadows labelled $1$, $2$, and $3$ in Figure
\ref{fig:m069}. The three $2$-handles of the Mom-$3$ structure meet
the cusp torus in a total of nine bridges, which correspond to the
highlighted edges in Figure \ref{fig:m069}. One can confirm easily
that the resulting lakes are simply connected, so this internal Mom-$3$
structure is full. Note that again the $1$-handles and $2$-handles of
this Mom structure have totally geodesic cores. It is also true (if
more difficult to confirm) that the submanifold corresponding to this
Mom-$3$ structure is the manifold $m202$, and that the components of the
complement of this submanifold consist of the original cusp
neighborhood and a single additional solid torus.

\section{Geometric and combinatorial Mom-$n$ structures}


As noted, the examples in the previous section have the property that
the $1$-handles and $2$-handles of $\Delta$ are all neighborhoods of
geodesic arcs and totally geodesic surfaces respectively. We will call
such an internal Mom-$n$ structure a \emph{geometric} Mom-$n$ structure
from now on to emphasize this fact. (Strictly speaking it would be more
correct to say ``geometric internal Mom-$n$ structure'' but for the sake
of concision we will assume that the word ``internal'' is implied.)
At the time of writing the authors know of no cusped manifold possessing an
internal Mom-$n$ structure which does not possess a geometric Mom-$n$
structure.


In addition to its geometric description, we can also describe a geometric Mom-$n$
structure in a one-cusped hyperbolic manifold $N$ combinatorially as follows.
Suppose $W$
is a choice of cusp neighborhood in $N$. Under the universal covering map
$\BH^3\rightarrow N$, the pre-image of $W$ is a collection of horoballs $\{B_i\}$,
and a $1$-handle in a geometric Mom-$n$ structure can be lifted to a
$1$-handle in $\BH^3$ joining two distinct horoballs $B_i$ and $B_j$.
Hence each $1$-handle in a geometric Mom-$n$ structure corresponds to an
element of the set of equivalence classes of unordered pairs
$\{(B_i,B_j)|i\not= j\}/\pi_1(N)$. Similarly each $2$-handle in a geometric
Mom-$n$ structure corresponds to an equivalence class of unordered triples
$(B_i,B_j,B_k)$ under the action of $\pi_(N)$, where $i$, $j$, and $k$ are
all distinct.

For each equivalence class in the set $\{(B_i,B_j)|i\not= j\}/\pi_1(N)$ the
\emph{ortho\-dis\-tance} $d(B_i,B_j)$ is well-defined; the set of all such
orthodistances counted with multiplicity forms the
\emph{orthodistance spectrum} $o(1)\le
o(2)\le\cdots$. By taking $T$ to be the boundary of a maximal embedded
cusp neighborhood we can assume that $o(1)=0$. We order the elements
of $\{(B_i,B_j)\}/\pi_1(N)$ by their orthodistances, and say that the
$n$-th \emph{orthopair class} $\mc{O}(n)$ is the equivalence class
with orthodistance $o(n)$. Each unordered triple of horoballs
$(B_i,B_j,B_k)$ is ``bounded'' by three orthopair classes
corresponding to the pairs $(B_i,B_j)$, $(B_j,B_k)$, and
$(B_k,B_i)$. We say a triple is a \emph{$(k,l,m)$-triple}, or is of
\emph{type $(k,l,m)$}, if the corresponding three orthopair classes
are $\mc{O}(k)$, $\mc{O}(l)$, and $\mc{O}(m)$ in some order.



With this, we can construct a combinatorial analogue to the idea
of a geometric Mom-$n$ structure:

\begin{definition}
A \emph{combinatorial Mom-$n$ structure} in a one-cusped hyperbolic
manifold $N$ is a collection of $n$ distinct equivalence classes of
triples of horoballs in the universal cover of $N$, of types
$(k_1,l_1,m_1)$, \ldots, $(k_n,l_n,m_n)$, such that the integers
$k_1$, $l_1$, $m_1$, \ldots, $k_n$, $l_n$, $m_n$, are all elements of
the same $n$-element subset of $\BZ_+$.
\label{defn:comb_Mom}
\end{definition}

We can say trivially that a manifold with a geometric Mom-$n$
structure must also possess a combinatorial Mom-$n$ structure. (For
example, the geometric Mom-$3$ structure on $m069$ described in the
previous section corresponds to a combinatorial Mom-$3$ structure with
triples of type $(1,1,2)$, $(1,3,3)$, and $(2,2,3)$.) Proving
the implication in the other direction is more difficult, as we will
see. But a motivation for doing so comes from the following:
\begin{theorem}
If $N$ is a one-cusped hyperbolic 3-manifold with $\Vol(N)\le 2.848$
then $N$ possesses a combinatorial Mom-$n$ structure with $n=2$ or
$3$.
\label{thrm:comb_Mom}
\end{theorem}

The complete proof of the above theorem is one of the key results of
\cite{GMM3}; we give a sketch of the result here. As with previous
results in this area the result comes from studying the maximal cusp
diagram of the manifold $N$.

Let $T$ be a cusp torus in $N$ such that the restriction of the metric
to $T$ is flat, and such that $T$ bounds a maximal cusp neighborhood
in $N$. Then $T$ will lift to a collection of horospheres in the
universal cover $\BH^3$, and the cusp neighborhood will lift to a
collection of horoballs $\{B_i\}$. Choose one of these horoballs and
call it $B_\infty$; then $B_\infty$ is fixed by a subgroup $H$ of
$\pi_1(N)$ isomorphic to $\BZ+\BZ$, and the quotient of $\partial
B_\infty$ by $H$ is just $T$. The cusp diagram of $N$ consists of $T$
together with the shadows of all of the other horoballs $\{B_i\}$ on
$\partial B_\infty$, modulo the action of $H$. These shadows form a
collection of circular disks on $T$ of varying radii, many of which
will overlap. We are interested in using these shadows to obtain a
lower bound on the area of $T$, for it is a simple matter to prove
that $\Vol(B_\infty/H)=\frac{1}{2}\Area(T)$, and $\mathrm{int}(B_\infty)/H$ is
embedded in $N$.

This result is used in \cite{Ada}, for example, to find the
minimum-volume cusped non-orientable manifold as follows. Consider all
of those horoballs which are tangent to $B_\infty$. For example, any
horoball $B_1$ such that $(B_1,B_\infty)$ is an element of $\mc{O}(1)$
would be such a horoball. Adams noted that there must be at least two
such horoballs modulo the action of $H$; call them $B_1$ and ${B_1}'$.
Specifically one can choose ${B_1}'$ to equal $g(B_\infty)$, where
$g\in\pi_1(N)$ is any element such that $g^{-1}(B_\infty)=B_1$. (The
fact that there is no element of $\pi_1(N)$ which exchanges $B_1$ and
$B_\infty$ implies that $B_1$ and ${B_1}'$ will not be in the same
$H$-orbit.) Since $B_1$ and ${B_1}'$ are disjoint horoballs both
tangent to $B_\infty$, their shadows on $\partial B_\infty$ must be
disjoint disks of radius $1/2$, which implies (using known bounds on
circle packing in the plane) that $\Area(T)\ge 2\pi(1/2)^2
(\sqrt{12}/\pi)=\sqrt{3}$. Hence from above we have
$\Vol(B_\infty/H)\ge \sqrt{3}/2$. Applying the bound on horoball
packing in $\BH^3$ due to Boroczky (\cite{Bor}) proves that $N$ must
have volume at least as large as that of a regular ideal
simplex. Showing that the Gieseking manifold is the unique
non-orientable manifold which achieves this minimum completes the
proof in \cite{Ada}.

Note that the above argument uses no facts at all about $N$ other than
the fact that it has a non-trivial cusp neighborhood.
By assuming orientability and carefully studying the
possible arrangements of horoball shadows on $T$, Cao and Meyerhoff
significantly extended Adams's argument, proving in \cite{CM} that the
figure-$8$ knot complement and its sibling were the smallest
orientable cusped hyperbolic $3$-manifolds. This result differs from
the result of \cite{Ada} in an interesting respect. In \cite{Ada},
Adams computed a lower bound which turned out to be realized by a
particular manifold (similarly, in \cite{Mey3} Meyerhoff computed
a lower bound on the volume of hyperbolic orbifolds that was realized
by a particular orbifold). In contrast, \cite{CM} sets up a
dichotomy. In \cite{CM}, maximal cusp diagrams are sorted according to
the following question: are there abutting full-sized disks (that is,
shadows of horoballs with $o(k)=0$) in the diagram? If not then there
is enough space to get good bounds on the area of the cusp torus and
hence the volume of the manifold. On the other hand, the presence of
such abutting disks is a special situation with group-theoretic
implications for the Kleinian group $\pi_1(N)$ that can be analyzed to
get specific worst cases, namely the two lowest-volume orientable
cusped hyperbolic $3$-manifolds.

There is a similar dichotomy at work in Mom technology. Loosely
speaking, in the maximal cusp diagram either the disks are not close,
which leads to good area and hence volume bounds, or else the disks
are close, leading to topological implications about the presence of
an internal Mom structure.  The key insight here is that the presence
or absence of a $(k,m,n)$-triple of horoballs in $N$ provides
geometric information about the arrangement of the corresponding
shadows. In particular suppose $(B_i,B_j,B_\infty)$ forms a triple of
type $(k,m,n)$, with the pair $(B_i,B_\infty)$ belonging to
$\mc{O}(k)$ and the pair $(B_j,B_\infty)$ belonging to
$\mc{O}(m)$. (Note that by transitivity of the group action, for any
triple of horoballs we can always assume that $B_\infty$ is one of the
elements of the triple and there are three ways to do so.) Then we can
prove (and do so in \cite{GMM3}) that the shadows of $B_i$ and $B_j$
on $\partial B_\infty$ have radii $\frac{1}{2}e^{-o(k)}$ and
$\frac{1}{2}e^{-o(m)}$ respectively, and that there is a path along
$T$ joining the centers of these two shadows of length
$e^{(o(n)-o(m)-o(k))/2}$. Let $e_n=e^{o(n)/2}$ for all $n$ for ease of
notation; then the shadows of $B_i$ and $B_j$ have radii
$\frac{1}{2}{e_i}^{-2}$ and $\frac{1}{2}{e_j}^{-2}$ respectively and
the path between their centers has length $e_n/(e_m e_k)$.

Now suppose in particular that $k=m=1$. Then $e_k=e_m=1$ since $T$ was
chosen to bound a maximal cusp neighborhood, so the shadows of $B_i$
and $B_j$ each have radius $\frac{1}{2}$, and the distance between
their centers is $e_n\ge 1$. This together with the obvious
circle-packing argument immediately gives a lower bound of $\sqrt{3}$
for the area of $T$, as in \cite{Ada}. However, we can say more.
We show in \cite{GMM3} that there are no $(1,1,1)$-triples of horoballs
(or indeed any $(m,m,m)$-triples for any $m$) in any orientable
manifold $N$, and hence $n\ge 2$ in the above computation. This
implies that the distance between the centers of these two shadows is
at least $e_2\ge 1$, and this in turn improves our bound on the area
of $T$ to $\sqrt{3}{e_2}^2$, a significant improvement if $e_2$ is
large. Whereas if $e_2$ is small, we can find at least two additional
shadows of radius $\frac{1}{2}{e_2}^{-2}$ which contribute area to $T$
and improve our lower bound in another way.

Suppose we make a further assumption, namely that there are no
$(1,1,2)$-triples in our manifold $N$. Then the distance between the
centers of the first two shadows is at least $e_3$, improving our area
bound still further. If we relax our assumption, and assume that there
is at most one $(1,1,2)$-triple up to the action of $\pi_1(N)$, then
the distance between the centers may be as low as $e_2$, but only
along at most one path in $T$. In other words, if we construct disks
of radius $e_3/2$ about the centers of the first two shadows, then
these new larger disks will overlap at most once and we can still
obtain a lower bound on the area of their union and hence on the area
of $T$. And if $N$ has two or more $(1,1,2)$-triples of horoballs,
then $N$ has a combinatorial Mom-$2$ structure.

Arguing in this fashion in \cite{GMM3} we show that a manifold
which does \emph{not} possess a combinatorial Mom-$2$ or Mom-$3$
structure must either have one of the parameters $e_2$ or $e_3$ be
sufficiently large (for example, if $e_3>1.5152$; the
significance of this number is explained in the next section) or else
it falls into one of $18$ cases enumerated by the presence or absense of
certain types of triples of horoballs. If $e_2$ or $e_3$ is sufficiently
large, an argument similar to the one above shows that $\Vol(M)>2.848$,
while in the remaining $18$ cases a rigorous computer-assisted version
of the above argument also shows that $\Vol(M)>2.848$, completing the
proof of Theorem \ref{thrm:comb_Mom}.

\section{Upgrading a combinatorial structure}

We now wish to take the data associated to a combinatorial Mom-$2$ or
Mom-$3$ structure in a manifold $N$ and use it to construct a
hyperbolic internal Mom structure. In short, we want to upgrade our
combinatorial structure to a topological one. The principle is
straightforward. If $(B_i,B_j,B_k)$ is a triple of horoballs in the
universal cover of $N$ which realizes a triple of type $(l,m,n)$ in a
combinatorial Mom-$2$ or Mom-$3$ structure, then there is a totally
geodesic $2$-cell $\sigma$ in $\BH^3$ bounded by the shortest geodesic arcs
from $B_i$ to $B_j$, $B_j$ to $B_k$, and $B_k$ to $B_i$, along with
arcs in the boundaries of all three horoballs. We wish to use the
projection of $\sigma$ to $N$ as the core of a $2$-handle in a
hyperbolic internal Mom-$2$ or Mom-$3$ structure for $N$. We do the
same for every triple in the combinatorial Mom-$n$ structure. Similarly
the geodesic arcs in the boundary of $\sigma$ project to
geodesic arcs in $N$, and we wish to use all such arcs as the cores of
$1$-handles in a hyperbolic internal Mom-$n$ structure.

There are numerous obstacles to this straightforward idea, however:
\begin{itemize}
\item The resulting handle structure may not be embedded in $N$.
\item Even if it is embedded, the boundary of the resulting
  submanifold $M$ may not be a collection of tori.
\item Even if $\partial M$ is a collection of tori, $M$ may not be a
  hyperbolic submanifold of $N$, i.e. $M$ may have an essential embedded
  sphere, annulus, or torus.
\end{itemize}
We now discuss how each of these obstacles can be overcome.


First there is the question of embeddedness. Here we make use of the
fact that the structure we are attempting to construct is a geometric
Mom-$n$ structure. For example, suppose a $1$-handle in our putative
Mom structure intersects itself, because the core geodesic arc of the
$1$-handle intersects itself. This implies that we have two horoballs
$B_i$ and $B_j$ in the universal cover such that the arc $\lambda$
from $B_i$ to $B_j$ intersects the arc $g(\lambda)$ for some
$g\in\pi_1(N)$. It can be shown that if $\lambda$ is sufficiently
short (specifically if $e_n\le \sqrt{2}$, where $\mc{O}(n)$ is the
corresponding orthopair class) then the interior of one of $B_i$ or $B_j$ must
intersect one of $g(B_i)$ or $g(B_j)$, a contradiction. For slightly
longer arcs (specifically if $e_n\le 1.5152$) it is shown in
\cite{GMM3} that this implies that the interior of one of $B_i$, $B_j$
intersects one of $g^k(B_i)$, $g^k(B_j)$ for some $k\le 4$. Note that
as mentioned previously it is also shown in
\cite{GMM3} that $\Vol(N)> 2.848$ if $e_2>1.5152$ or $e_3>1.5152$,
which implies that longer arcs need not be considered.

Another possibility is that two different $1$-handles intersect, in
which case there exists four horoballs $B_i$, $B_j$, $B_k$, and
$B_l$ in the universal cover such that the arc from $B_i$ to $B_j$
intersects the arc from $B_k$ to $B_l$. In this case it is shown in
\cite{GMM3} that if the arcs are
sufficiently short (again meaning that $e_m\le 1.5152$ and $e_n\le
1.5152$ where $\mc{O}(n)$ and $\mc{O}(m)$ are the appropriate
orthopair classes) then a new, simpler combinatorial Mom structure can be
constructed which excludes one of the two troublesome $1$-handles.

Other types of intersections may occur between the various handles of
our putative handle structure, but in each case one of these three things
occurs: either we obtain a geometric contradiction, or the
lengths of the $1$-handles in the structure assume values which imply
that $\Vol(N)>2.848$, or else we can construct a strictly simpler
combinatorial Mom structure and start again. Therefore by induction we
have the following:
\begin{theorem}
If $N$ is a one-cusped hyperbolic 3-manifold and $\Vol(N)\le 2.848$
then $N$ has a combinatorial Mom-$n$ structure with $n=2$ or $3$
corresponding to a handle structure on an embedded submanifold $M$.
\end{theorem}
The arguments appear in full detail in \cite{GMM3}; they are lengthy
but elementary.


The next concern is the topology of the boundary components of the
resulting manifold $M$. This turns out to be the simplest problem to
overcome. By construction the Euler characteristic of $M$ is $0$, and
hence if the boundary components are not all tori then one of them
must be a sphere. Here again we can take advantage of the fact that
our putative Mom structure is a geometric one. Given a handle
structure of the type described up to this point, with only two
$2$-handles whose cores are totally geodesic $2$-cells, it is in fact
impossible to construct a submanifold $M$ of $N$ with at least two
boundary components, one of them a sphere. Hence if we start with a
combinatorial Mom-$2$ structure then the boundary components of $M$
are automatically tori. With a combinatorial Mom-$3$ structure there
is one way to construct a submanifold $M$ with a spherical boundary
component, but that way requires that the combinatorial structure have
exactly two triples of type $(k,l,m)$ where $k$, $l$, and $m$ are
three distinct integers. We say that such a combinatorial structure
is not \emph{torus friendly}.
The same analysis described in the previous
section shows that a manifold with such a combinatorial Mom-$3$
structure, and no torus friendly combinatorial Mom structures satisfying
the previous constraints on $e_2$ and $e_3$,
must have volume greater than $2.848$, thus dealing with this one
exceptional case.


Now we have a submanifold $M$ embedded in $N$, such that the boundary
of $M$ is a union of tori, together with a handle structure of the
appropriate type. If we can show that $M$ can be assumed to be
hyperbolic then we will have the hyperbolic internal Mom-$n$ structure
that we desire. (Astute readers may have noticed that
we haven't shown that each
component of $\partial M$ is either the cusp torus $T$ or else
bounds a solid torus in $N\setminus i(M)$. However if $M$ is
hyperbolic then this condition will be satisfied automatically: given
the way that $M$ was constructed, the only other possibility is that
some boundary component of $M$ bounds a ``tube with knotted hole'' in
$N$, and if this is the case then $M$ would contain an embedded
essential sphere.)

As with embeddedness, the idea is to show that if $M$ is not
hyperbolic then we can find a strictly simpler internal Mom structure, and
hence we can assume that $M$ is hyperbolic by induction. For example,
suppose $M$ contains an embedded essential sphere. Following the ideas
of Matveev, we assume the sphere is in normal position with respect to
the handle structure of $M$ and then split both $M$ and the handle
structure along the surface. After throwing away the component of the
split manifold which does not contain $T$, the result is a submanifold
of $N$ with torus boundary components and one spherical boundary
component, which must bound a ball in $N$.  Adding that ball as a
$3$-handle and cancelling it with a $2$-handle results in a new
submanifold $M_1$ with torus boundary components, and a new handle
structure $\Delta_1$ which has strictly lower complexity (in the sense
of Matveev) than the handle structure we started with. The new handle
structure is not necessarily in the form of a Mom-$n$,
because the $2$-handles of the structure may not be attached to the
correct number of $1$-handles. But $2$-handles which attach to four or
more $1$-handles can be split into $2$-handles of valence three. And
$2$-handles which attach to two or fewer $1$-handles can be eliminated
from the handle structure, usually by cancelling them with a
$1$-handle. Furthermore, these operations do not increase Matveev
complexity. The end result is a handle structure for $M_1$ in the form
of a Mom-$k$, resulting in an internal Mom-$k$ structure on $N$
where the new Mom number $k$ is strictly less than the Mom number we
started with.

Since $M$ is compact and $\partial M$ is a non-empty union of tori,
$M$ is either hyperbolic or else contains an essential embedded 2-sphere,
annulus, or torus.  In each case we can construct a new manifold
$M_1\subset N$ with a simpler handle structure as measured by Matveev
complexity; the spherical case is outlined in the previous paragraph, 
but the other cases are more complicated. A technical point is that the
new handle structure might not have simply connected lakes
(i.e. is not full) and hence this proof requires $n\le 4$ since it
relies on a solution to Problem \ref{general based} discussed in \S 10 which is
only known for $m\le 2$.  See \cite{GMM2} for more details.

Hence we conclude:
\begin{theorem}
If $N$ is a one-cusped hyperbolic 3-manifold with $\Vol(N)\le 2.848$
then $N$ has an internal Mom-$n$ structure $(M,T,\Delta)$ with $n=2$
or $3$ such that $M$ is hyperbolic, and hence $N$ is obtainable by a
Dehn filling on $M$.
\label{thrm:vol_to_Mom}
\end{theorem}

\section{Classification of Mom-$n$'s, $n < 4$}

Theorem \ref{thrm:vol_to_Mom} allows us to enumerate all one-cusped
manifolds with $\Vol(N)\le 2.848$ in two steps. First, we need to
enumerate all hyperbolic Mom-$2$'s and Mom-$3$'s, and second, we need
to enumerate all Dehn fillings on those manifolds which produce
one-cusped manifolds with volume less than or equal $2.848$.

The first step is the more straightforward one. If $(M,T,\Delta)$ is a
hyperbolic Mom-$n$, then in particular it must be full. This implies
that $\Delta$ retracts to a \emph{spine} for
$M$, i.e. a cellular complex which intersects every homotopically
non-trivial simple closed curve in $M$. This spine
consists of:
\begin{itemize}
\item a $0$-cell corresponding to every island,
\item a $1$-cell corresponding to the core of every $1$-handle of
  $\Delta$, and a $1$-cell corresponding to every bridge, and
\item a $2$-cell corresponding to the core of every $2$-handle of
  $\Delta$, and a $2$-cell corresponding to every lake.
\end{itemize}
Note that if $(M,T,\Delta)$ were not full, then there might be lakes
which are not disks and this construction would not work.  We will
abuse notation and let $\Delta$ refer to the spine of $M$ from
this point forward. 

Being a spine of $M$, $\Delta$ is dual to an ideal polyhedral cellulation
of $M$ with one $3$-cell for every $0$-cell of $\Delta$ and an ideal
vertex for every boundary component of $M$; clearly there
are $2n$ such $3$-cells, two for each $1$-handle in the Mom-$n$. If
a $1$-handle has valence two, then the corresponding dual $3$-cells
will be pyramids built on a digonal base; we ignore such $3$-cells after
flattening them down to a face. If a $1$-handle has valence $v>2$,
then the two corresponding dual $3$-cells are each pyramids built on a
$v$-sided polygon as base; see Figure
\begin{figure}
\begin{center}
\includegraphics{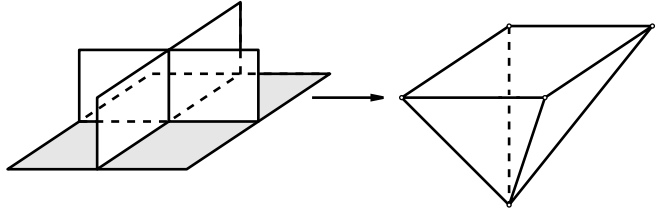}
\end{center}
\caption{A neighborhood of an endpoint of a $1$-handle of valence $v$
together with the dual $v$-sided ideal pyramid, where $v=4$.}
\label{fig:dual_cells}
\end{figure}
\ref{fig:dual_cells}. Furthermore, the
bases are each dual to the $1$-handle itself, and hence are glued
together in the resulting dual cellulation, resulting in a $v$-sided
dipyramid. Also, the sums of the valences of the $1$-handles in the
Mom-$n$ equals the sums of the valences of the $2$-handles,
which is $3n$ by definition.

This generates a finite number of possible polyhedral decompositions
for a hyperbolic Mom-$2$ or Mom-$3$. Specifically, a hyperbolic
Mom-$2$ must be composed of either a single ideal octahedron, or two
ideal three-sided dipyramids. A hyperbolic Mom-$3$ must consist of
either a single ideal five-sided dipyramid, an ideal octahedron and an
ideal three-sided dipyramid, or three ideal three-sided
dipyramids. And furthermore, the ``polar'' ideal vertices of these
dipyramids, i.e. the vertices not adjacent to the bases of the
original pyramids, are all dual to the $T\times \{0\}$ boundary
component of $M$ and
hence must be identified to other ``polar'' vertices by the
face-pairing identifications. This is enough information to enumerate
all of the possible hyperbolic Mom-$2$'s and Mom-$3$'s.

Computer assistance was used here, to generate a candidate list of
polyhedral gluings. The number of such gluings is a factorial function
of the total number of faces; to shrink this list, at this stage the
computer was programmed to check that the links of ideal vertices had
Euler characteristic zero, and to use the obvious symmetries of the
dipyramids to eliminate redundant gluings whenever possible. This
resulted in $44$ candidate polyhedral gluing descriptions for
Mom-$2$'s and $4187$ such descriptions for Mom-$3$'s. Then SnapPea was
used to make a preliminary identification of which gluing descriptions
corresponded to hyperbolic manifolds, and to eliminate duplicates
among the hyperbolic manifolds.

In each case where SnapPea found a hyperbolic structure it also found
a triangulation of the manifold identical to a triangulation of a
manifold in the SnapPea census. Since Harriet Moser has confirmed the
hyperbolicity of all the manifolds in the census (\cite{Mos}) this
confirms the hyperbolicity of our Mom-$2$'s and Mom-$3$'s; i.e.,
SnapPea did not make any false positive errors. For those manifolds
for which SnapPea failed to find a hyperbolic structure, the
fundamental group of the manifold was computed and examined by hand,
and it was shown that none of these groups was the fundamental group
of a finite-volume hyperbolic manifold. (In most cases, this was accomplished
by showing that the group had a non-trivial center. A few cases required
a more detailed analysis, usually involving finding distinct maximal abelian
subgroups with non-trivial intersection, which is impossible in a discrete
co-finite volume group of hyperbolic isometries.)
In this way, we confirmed that
SnapPea did not make any false negative errors. The results of this
analysis are as follows:

\begin{theorem}
If $(M,T,\Delta)$ is a hyperbolic Mom-$2$ then $M$ is homeomorphic to
one of the manifolds known in the SnapPea census as $m125$, $m129$, or
$m203$. If $(M,T,\Delta)$ is a hyperbolic Mom-$3$ then $M$ is
homeomorphic to one of the preceeding three manifolds, or to one of $m202$,
$m292$, $m295$, $m328$, $m329$, $m359$, $m366$, $m367$, $m391$,
$m412$, $s596$, $s647$, $s774$, $s776$, $s780$, $s785$, $s898$, or $s959$.
\label{thrm:Mom_enum}
\end{theorem}\qed

Some comments about the above list: $m129$ is the complement of the
Whitehead link, and $m203$ is the complement of the link known as
$6_2^2$ in standard link tables. Among the Mom-$3$'s, $s776$ is the
complement of the alternating three-element chain link, sometimes
known as the ``magic manifold'' and extensively analyzed in \cite{MP};
it is also the only three-cusped manifold in the above list, and the
manifolds from $m125$ to $m391$ in the above list can all be obtained by a
Dehn filling on $s776$.

A similar analysis was performed for hyperbolic Mom-$4$'s, and
resulted in a list of 138 hyperbolic Mom-$4$'s which included each of
the 21 manifolds listed above. Some aspects of this analysis, however,
are still conjectural. In particular we have not yet analyzed the
fundamental groups of all the manifolds for which SnapPea failed to
find a hyperbolic structure. The list of Mom-$4$'s includes some
manifolds which do not appear in the SnapPea census; descriptions of
these manifolds can be found in \cite{GMM2}. See the Mom-$4$ Conjecture
in \S 10.

\section{Enumeration of filled manifolds}

After identifying all hyperbolic Mom-$2$'s and Mom-$3$'s, the
next step is to identify all Dehn fillings on those manifolds which
might result in one-cusped hyperbolic manifolds of volume less than or equal to
$2.848$. We turn to \cite{FKP}, which says the following:

\begin{theorem}(Futer, Kalfagianni, and Purcell) Let $M$ be a
  complete, finite-volume hyperbolic 3-manifold with cusps. Suppose $C_1$,
\ldots, $C_k$ are disjoint horoball neighborhoods of some subset of
the cusps. Let $s_1$, \ldots, $s_k$ be slopes on $\partial C_1$,
\ldots, $\partial C_k$, each with length greater than $2 \pi$. Denote
the minimal slope length by $l_{\mr{min}}$. If $M(s_1,\ldots,s_k)$
satisfies the Geometrization Conjecture, then it is a hyperbolic
manifold, and
\[
\Vol(M(s_1,\ldots,s_k)) \ge \left(1-\left(\frac{2\pi}
{l_{\mr{min}}}\right)^2\right)^{3/2} \Vol(M).
\]
\label{thrm:fkp}
\end{theorem}

Note that a \emph{slope} here just denotes the homotopy class of curves on
the boundary of $M$ which end up glued to the boundary of a disk after
Dehn filling, and $M(s_1,\ldots,s_k)$ denotes the filled manifold.
Assuming the Geometrization Conjecture is true, then the above theorem
has the following simple reformulation, as noted in \cite{M}:

\begin{corollary}Suppose $M$, $s_1$, \ldots, $s_k$, and $l_{\mr{min}}$
are defined as above, and let $N=M(s_1,\ldots,s_k)$. If $N$ is
hyperbolic we have
\[
l_{\mr{min}} \le
2\pi \left( \sqrt{1- \left(\frac{\Vol(N)}{\Vol(M)}\right)^{2/3}}
\right) ^{-1}\,.
\]
\label{cor:fkp}
\end{corollary}

(Note that the right-hand side of the above inequality is always
greater than or equal to $2\pi$, so that condition on the length of
the slopes can dropped; if the boundary slopes do not all have length
greater than $2\pi$ then the inequality is trivially true.)

Applying Corollary \ref{cor:fkp} to each of the manifolds in Theorem
\ref{thrm:Mom_enum}, and assuming that $\Vol(N)\le 2.848$, we see that
there are a finite number of Dehn fillings on each manifold that need
to be considered. Furthermore these fillings are easily
enumerated; see \cite{M} for details. It is then necessary to
determine which of these fillings result in hyperbolic manifolds and
which do not, and for those manifolds which are hyperbolic it is
necessary to rigorously prove which filled manifolds actually have
volume no greater than $2.848$. Again, these steps were completed with
computer assistance; first, Snap (see \cite{Go}) was used to make a preliminary
determination, and then Snap's conclusions were confirmed by more
rigorous means.

We used Snap rather than SnapPea at this step to make more extensive
use of Harriet Moser's algorithm in \cite{Mos}. That is, rather than
attempt to locate each hyperbolic manifold in SnapPea's census (a more
daunting task here than in the proof of Theorem \ref{thrm:Mom_enum}
due to the larger number of cases to consider) we instead applied
Moser's algorithm directly to the filled manifolds under consideration
to prove their hyperbolicity. Moser's algorithm is designed to use
data produced by Snap as input, hence it was more convenient here to
use Snap rather than SnapPea. In this way we confirmed the
hyperbolicity of those one-cusped filled manifolds for which Snap
claimed to find a hyperbolic structure.

For those manifolds for which Snap failed to find a hyperbolic
structure, we confirmed non-hyperbolicity by examining the manifolds'
fundamental groups and checking for non-trivial centers and the like,
as in the proof of Theorem \ref{thrm:Mom_enum}.

Next, for those filled manifolds which are hyperbolic it is
necessary to rigorously establish which manifolds satisfy $\Vol(N)\le
2.848$. Here again we turn to computer assistance. At this point in
the process Snap has already produced an ideal triangulation of all of
the manifolds in question, and has computed values for the shape
parameters of each of the tetrahedra involved. From this information,
the volume of each manifold can be computed via the Lobachevsky
function $\text{\cy{L}}(\theta)$ (see, for example, \cite{Mil}); however the question of
floating-point error must be addressed. Fortunately, one of the
intermediate steps in Moser's algorithm is to compute an error bound
$\delta$ on the distance in $\BC^k$ between the shape parameters
computed by Snap and the actual shape parameters associated to the
hyperbolic metric. With this information and using affine 1-jets as in
\cite{GMM3} we can rigorously compute an upper and lower bound on the
volume of each hyperbolic manifold under consideration.

Finally, for each hyperbolic filled manifold of sufficiently small
volume an isometry was found between that manifold and a manifold in
the SnapPea census. By this process, the following theorem was proved
in \cite{M}:

\begin{theorem}
If $N$ is a one-cusped hyperbolic 3-manifold with $\Vol(N)\le 2.848$
and $N$ is obtainable by filling one of the manifolds in Theorem
\ref{thrm:Mom_enum}, then $N$ is one of the first ten orientable
one-cusped manifolds in the SnapPea census; that is, $N$ is one of
$m003$, $m004$, $m006$, $m007$, $m009$, $m010$, $m011$, $m015$, $m016$,
or $m017$.
\label{thrm:small_fillings}
\end{theorem}\qed

This together with Theorem \ref{thrm:vol_to_Mom} establishes the
following:

\begin{theorem}
There are only ten orientable one-cusped hyperbolic manifolds $N$ with
$\Vol(N)\le 2.848$, and these are exactly the first ten orientable
one-cusped manifolds in the SnapPea census.
\label{thrm:small_cusped_vols}
\end{theorem}\qed

\section{History of the compact case}

Before completing our discussion of the proof of Theorem \ref{main},
we take a moment to discuss the history of the search for the
minimum-volume compact hyperbolic $3$-manifold.

The Kazhdan-Margulis theorem when applied to
hyperbolic 3-manifolds
establishes the existence of a constant $\epsilon$ and a natural
decomposition of any hyperbolic $3$-manifold into an $\epsilon$-thick
part and an $\epsilon$-thin part; see \cite{KaMa}. This decomposition can be used to
show that there is a positive lower bound to the volume of hyperbolic
$3$-manifolds. This approach is carried out in Section 4 of \cite{Bu},
where Buser and Karcher implemented an idea of Gromov.

Virtually all low-volume bounds arise out of analyses of embedded
solid tubes. The first such analysis was in \cite{Mey1}, where
Meyerhoff used Jorgensen's
trace inequality (see \cite{Jo}) to construct solid tubes around short
geodesics in hyperbolic $3$-manifolds. The shorter the geodesic, the
larger the solid tube is, not only in radius but in volume.  From this
construction follows a trade-off argument: either a hyperbolic
3-manifold has a short geodesic and hence a solid tube with decent
volume, or it doesn't, in which case there must be an embedded ball,
again of decent volume.  The resulting bound is roughly $0.0006$.

Subsequently, sphere-packing in $\BH^3$ was used to gain some control
over the volume of a hyperbolic $3$-manifold outside of an embedded
ball.  This improved volume bound in the embedded ball case can be
used to make an improved trade-off in the tube-versus-ball argument of
\cite{Mey1}.  The necessary sphere-packing results had been produced
earlier by Boroczky and Florian (see \cite{BF}, or \cite{Bor} for an
English version), and the volume bound was pushed up to roughly
$0.0008$ in \cite{Mey2}.  Gehring and Martin noted that the
calculations in \cite{Mey1} could be made a bit finer and produced a
bound of roughly $0.001$ in \cite{GM1}.

For a number of years no further progress on lower bounds was made in
the closed case.  However in the cusped case, Adams was able to use
the maximal cusp diagram to get improved bounds for $v_\omega$. As
previously noted in \S 4, Adams showed in \cite{Ada} that the volumes
of orientable \emph{and} non-orientable cusped hyperbolic
$3$-manifolds are bounded below by $\mc{V}=1.01\ldots$, the volume of
the ideal regular simplex; he further showed that in the
non-orientable case this volume is realized uniquely by the Gieseking
manifold.

In the mid-1990's Gabai, Meyerhoff and N. Thurston needed to greatly
improve known solid-tube radius bounds so as to apply Gabai's Rigidity
Theorem and prove that homotopy hyperbolic $3$-manifolds are
hyperbolic. To do this they analyzed $2$-generator groups naturally
associated to the shortest geodesic in a hyperbolic $3$-manifold. The
space of such groups is determined by three complex parameters, two of
which are the complex length of the shortest geodesic and the complex
distance between the two nearest lifts of that geodesic to $\BH^3$.
(A third parameter is necessary to actually determine the isometry
between those two lifts.)
With some work, the relevant $3$-parameter space was shown to
be compact and a rigorous computer analysis of this space was carried
out in \cite{GMT}.  One crucial tool needed to obtain compactness of
the parameter space is the solid-tube construction described above; in
particular, according to \cite{Mey1}, shortest geodesics of length
less than $0.0979$ automatically have big enough tubes.  The result of
the \cite{GMT} analysis shows that, with seven manageable families of
exceptions, the shortest geodesic in a hyperbolic $3$-manifold has
tube radius greater than $\log(3)/2$. This is a considerable
improvement on previous tube bounds, and directly leads to a manifold volume
bound of roughly $0.1$, a hundred-fold improvement (the exceptions all
must have volume above $1.0$).

Coupling this $\log(3)/2$ result with a (slightly earlier) remarkable theorem of Gehring and Martin
produced a volume bound of $0.166$.  See \cite{GM2} and \cite{GMT}.  Gehring and Martin
showed how to generalize the maximal cusp diagram approach (e.g., in
\cite{Ada}) to the case of closed manifolds and maximal solid
tubes. In the cusped case the maximal cusp diagram is formed by the
shadows of horoballs $\{B_i\}$ on the boundary of a base horoball
$B_\infty$; the resulting shadows are disks. In the closed case,
instead of horoballs a maximal solid tube around a short geodesic is
lifted to $\BH^3$ to get a collection of solid tubes. One of these
tubes is designated as the base, and then the other tubes are
projected to the boundary of the base tube thereby creating a
collection of shadows on that boundary.  Arguments similar to those in
the cusped case should work here as well, but there are daunting
technical problems involved in projecting tubes onto the tube
boundary, as the shapes of the shadows are now much more complicated.
Gehring and Martin were able to sidestep many of these problems by
embedding a ball within each tube and projecting that ball to the base
tube instead, although this sacrifices some volume. Their approach
results in a formula for volume of the solid tube in terms of the
radius of the tube alone; surprisingly the length of the core geodesic
is not needed.

Przeworski then improved on the Gehring-Martin maximal tube diagram
approach by analyzing the shadows of the solid tubes themselves rather
than the shadows of balls, and by analyzing some of the volume outside
the maximal tube. He 
produced a volume lower bound of $0.28$, in \cite{Prz0} and
\cite{Prz1}.

Given the success of the \cite{GMT} method, it seemed natural to try
to extend the parameter space argument of \cite{GMT} and to focus it
more strictly on volume questions.  The first step in such a procedure
would be to produce an appropriate compact parameter space to analyze.
Because the lowest known volume manifold was the Weeks manifold with
volume $0.942\ldots$ and because the volume bound produced by Adams in
the cusped case is $1.01\ldots$ it would seem clear that very short
geodesics (nearly cusps) could be eliminated from the parameter space
argument.  The problem was that Adams's bound utilized
horoball-packing results, and these packing results could not be
generalized at the time to the case of short geodesics and big tubes.
Thus the volume bound of \cite{GM2} in the case of very short
geodesics (or tubes of very large radius) approaches $\sqrt{3}/2 =
0.866\ldots$ in the limit. This is Adams's result when horoball-packing
is ignored, and not adequate for attaining the low-volume manifold.

The authors were able to produce the desired compact parameter space
without use of tube-packing in \cite{GMM1}.  We introduced a simple
method for improving on the cusped volume bound of Adams, and then
perturbed the method to the closed case thereby producing a compact
length bound for the parameter space.  Our method was to look at the
next-largest disks in the maximal cusp diagram, i.e. the shadows of
the horoballs which are at a distance of $o(2)$ from $B_\infty$, to use the
language of this paper. What makes this argument work is the following
dichotomy, similar to the arguments of \S 4: either these horoballs
are close to $B_\infty$ and the associated shadows produce
substantial extra area and hence volume, or they are not close in
which case the centers of the full-sized disk shadows must be far from
each other, again producing substantial extra area and hence volume.
This approach carries over to the closed case as well.

Note that despite having an appropriate compact parameter space to
analyze, the \cite{GMT} approach has so far proven difficult to extend
beyond the original bounds.  Note also that the implicit improvement
on Adams's cusp bound in \cite{GMM1} was considerably less than the
improvement that had already been attained by Cao and Meyerhoff.  In
fact, Cao and Meyerhoff had doubled Adams's bound and this turned out to be
precisely $v_\omega$; see \cite{CM}.  One could try to perturb the
Cao-Meyerhoff methods to the closed case, but they are quite intricate
and this approach was sidelined as other events moved to the fore;
specifically, \cite{GMM1} seemed to spark a flurry of activity.

First Marshall and Martin showed how to rethink certain aspects of the
Gabai-Meyer\-hoff-Milley argument and were able to improve the volume
bound for closed manifolds to $0.2855$, in \cite{MM1}.  In separate
research, Marshall and Martin developed tube-packing methods in
$\BH^3$ in the large-tube-radius setting.  That is, if the radius of a
tube is sufficiently large (roughly radius $5$) in a hyperbolic
$3$-manifold, then the lifts of the tube are geometrically
sufficiently similar to a horoball that packing bounds can be
successfully obtained; see \cite{MM2}.

Then by an elegant argument in \cite{Prz2}, Przeworski obtained
tube-packing results in broader generality than \cite{MM2}.
Przeworski's results are typically applied whenever any new volume
result is established, via a tube volume argument.

In a major development, Agol studied the relationship between a closed
hyperbolic $3$-manifold with an embedded geodesic and the associated
cusped manifold obtained by removing that geodesic.  Using delicate
geometric constructions and applying a result of Boland-Connell-Souto \cite{BCS},
Agol was able to bound the volume of the closed manifold in terms of
the volume of the cusped manifold and the radius of the maximal tube
around the geodesic.  Using the volume bounds of \cite{CM}, Agol
produced a volume lower bound of $0.32$ for closed manifolds in
\cite{A1}. Przeworski then followed this result with \cite{Prz3},
which further improved bounds for the density of tube
packings. Combining these results with the results in \cite{A1},
Przeworski obtained a volume bound of $0.3324$.

Finally, Agol and Dunfield realized that Perelman's work on Ricci
curvature (en route to Perelman's proof of Thurston's Geometrization
Conjecture) substantially improves the results of \cite{A1}, which
involved Ricci curvature arguments.  The volume bound so produced is
$0.66$, and Agol, Storm, and Thurston re-proved this result in
\cite{AST}.  (Both the $0.3324$ and $0.66$ results use the Log(3)/2 theorem of \cite{GMT}.)  Of course, this is close to the hoped-for bound of
$0.942...$. 

This is the point where Mom technology re-enters the narrative,
for it is Agol and Dunfield's result together with the classification
of low-volume cusped manifolds due to Mom technology that allow us to
prove that the Weeks manifold is volume-minimizing, as described in
the next section.

\section{The minimum-volume closed manifold}

Starting with Theorem \ref{thrm:small_cusped_vols}, we want to apply
the results of Agol and Dunfield to tackle the closed case. Although
these results are documented in \cite{AST}, for our purposes it was
convenient to use a slightly different formulation of the same result
which appears as Lemma 3.1 of \cite{ACS1}:

\begin{lemma}
Suppose that $N$ is a closed orientable hyperbolic $3$-manifold and
that $C$ is a shortest geodesic in $N$ with an embedded tubular
neighborhood of radius at least $\log(3)/2$. Set $M=N\setminus C$, equipped
with a hyperbolic metric. Then $\Vol(M)<3.02 \Vol(N)$.
\label{lem:acs}
\end{lemma}

Recall that the Weeks manifold $W$ satisfies $\Vol(W)=0.9427\ldots$,
and furthermore note that in the above theorem if $\Vol(M)>2.848$ then
$\Vol(N)>0.943>\Vol(W)$. This is the reason behind the choice of
$2.848$ as a volume bound in Theorem
\ref{thrm:small_cusped_vols}. Also note that if the shortest geodesic
in $N$ does not have an embedded tubular neighborhood of radius at
least $\log(3)/2$ then $\Vol(N)>\Vol(W)$ according to
\cite{GMT}. Hence combining Theorem \ref{thrm:small_cusped_vols},
Lemma \ref{lem:acs}, and \cite{GMT} yields the following:

\begin{theorem}
Suppose that $N$ is a closed orientable hyperbolic $3$-manifold with
$\Vol(N)\le \Vol(W)$. Then $N$ is obtained by Dehn filling on one of
the first ten orientable one-cusped manifolds in the SnapPea census.
\label{thrm:small_closed_fillings}
\end{theorem}\qed

Clearly we now need to enumerate all Dehn fillings of those ten
manifolds which can result in a closed manifold with volume no greater
than $0.943$. This analysis, performed in \cite{M}, is similar enough
to the proof of Theorem \ref{thrm:small_fillings} that we will not
repeat the details here. We will mention a complication which did
not occur in the proof of Theorem \ref{thrm:small_fillings},
however. One of the closed manifolds that needs to be examined for
this analysis is the manifold $\Vol 3$, the third-smallest known closed
hyperbolic manifold. As the name suggests, $\Vol 3$ does not have volume
smaller than the Weeks manifold. However proving this using the
techniques used in the proof of Theorem \ref{thrm:small_fillings} is
complicated by the fact that $\Vol 3$ is
the only known hyperbolic manifold for which a non-negatively oriented ideal
triangulation has not been found. Moser's algorithm and the standard
formula for hyperbolic volume of manifolds both depend on having an
ideal triangulation without negatively oriented tetrahedra. Hence in
this one case we considered not $\Vol 3$ but the unique double
cover of $\Vol 3$, and showed that its double cover had volume no less
than $1.886$. (See also the discussion following Problem \ref{prob:10_34}
in the next section.)
All other closed fillings of the manifolds listed in
Theorem \ref{thrm:small_fillings}, except for the Weeks manifold, have
volume greater than $0.943$. This completes the proof of Theorem \ref{main}.

\section{Problems and Directions}

We close this paper with some open problems and possible future
directions for research. To begin with, the authors view the work of
\cite{GMM1}, \cite{GMM2}, \cite{GMM3} and \cite{M} as steps in
addressing the following:

\begin{ComplexityConjecture} (Thurston, Hodgson-Weeks, and
Matveev-Fom\-enko) The complete low-volume hyperbolic $3$-manifolds
can be obtained by filling cusped hyperbolic $3$-manifolds of small
topological complexity.\end{ComplexityConjecture}

\begin{remark} A detailed discussion of this conjecture can be found
in the introduction to \cite{GMM2}.  In particular, one of the
challenges is to quantify the adjectives \emph{low} and \emph{small}.
Our point of view is that, at least for low-volume manifolds, the Mom
number is an excellent topological measure compatible with
volume.

The experimental evidence provided by SnapPea is compelling.  Among
the $1$-cusped manifolds in the SnapPea census, experimental evidence
suggests that all such manifolds with volume at most $3.18$ (resp. $4.05$,
resp. $5.33$) have Mom number two (resp. at most three, resp. at most
four).

Among the $117$ smallest closed orientable $3$-manifolds in the
census, i.e. the manifolds with volume less than $2.5$, all but $5$
have internal Mom-$2$ structures which are based on a shortest closed
geodesic.  This means that some boundary component of the Mom manifold
bounds a solid torus whose core is a shortest geodesic. Four of the
remaining $5$ have internal Mom-$3$ structures based on a shortest
geodesic; they are m$038(-1,2)$, m$038(1,2)$, m$038(4,1)$, and m$038(-5,1)$.
The remaining manifold, m$207(1,2)$ in the SnapPea census, has volume
approximately equal to $2.468$ and an internal Mom-$4$ structure based
on a shortest geodesic (specifically, the core geodesic is a shortest
geodesic and m$207$ itself has a Mom-$4$ structure).
\end{remark}

\begin{problem} Develop a Mom technology theory for closed orientable
hyperbolic $3$-manifolds and directly prove that the Weeks manifold is
the lowest volume closed hyperbolic $3$-manifold.\end{problem}

That is, generalize the methods of \cite{GMM3} to directly address
closed $3$-manifolds.  For example, let $\gamma$ be a shortest
geodesic in the $3$-manifold $N$ and $V$ a maximal solid tube about
$\gamma$.  By passing to the universal covering of $M$, fixing one
preimage $V_0$ of $V$ and considering the other preimages $\{V_i\}$ it
makes sense to talk about orthoclasses, triples and hence the notion
of combinatorial Mom-$n$ structure. As noted in \S 8, Prezworski has
\cite{Prz1} developed the theory of shadows of projections of one
solid tube onto another in \cite{Prz1}.  With a generalization of the
``lessvol'' function which appears in \cite{CM} and \cite{GMM3} to the
closed case, one could directly generalize \cite{GMM3}.

\bigskip
If $V$ is the maximal solid tube described above, by \cite{GMT} we
know that either $V$ has tube radius at least $\log(3)/2$ or $M$ lies
in one of seven exceptional families of $3$-manifolds.  By \cite{JR}
two of these families are isomorphic, and by \cite{GMT}, \cite{Ly} and
\cite{CLLMR} associated to each of these families is a unique closed
orientable hyperbolic $3$-manifold.

\begin{problem} Complete the proof of Conjecture 1.31 of \cite{GMT}
by showing that each of these six manifolds $N_0,
N_1,\cdots, N_5$ cover only themselves.\end{problem}

\begin{remark} Jones and Reid showed in \cite{JR} that $N_3$ (also known
as $\Vol 3$) nontrivially covers no $3$-manifold and Reid showed in an
appendix to \cite{CLLMR}
that $N_1$ and $N_5$ nontrivially cover no $3$-manifold.

See Corollary 1.29, Remarks 1.32, and Theorem 4.1 of \cite{GMT} for more
information about these exceptional manifolds.  \end{remark}

\begin{problem} Find the lowest volume closed non\-orientable
3-manifolds.\end{problem}

One difficulty with nonorientable $3$-manifolds is that one of the
hypotheses of the $\log(3)/2$ theorem of \cite{GMT} is
orientability. In particular it is not applicable to orientation
reversing curves.  Nevertheless, by Milley some information carries
over; see \cite{M1}.

\begin{problem} Find a version of the $\log(3)/2$ theorem for shortest
geodesics in closed non\-orientable hyperbolic $3$-manifolds and for
cusped hyperbolic $3$-manifolds.  Improve the value of $\log(3)/2$ for
closed orientable $3$-manifolds.\end{problem}

\begin{problem} Find the lowest volume cusped non\-orientable
$3$-manifolds with torus cusps.\end{problem}

\begin{remark} \cite{Ada} showed that the Gieseking manifold is
the unique lowest volume cusped, non\-orientable hyperbolic
$3$-manifold, but this manifold has a Klein bottle cusp. The manifold
known as $m131$ in the SnapPea census is a non\-orientable manifold with
a single torus cusp and has volume equal to the volume of the
Whitehead link complement, which is approximately $3.663$; it is
the smallest such manifold known to the authors. The filled 
non\-orientable manifold $m131(3,1)$ has volume equal to the volume
of the figure-eight knot complement, i.e. approximately $2.029$, and
is the smallest compact non-orientable manifold known to the authors.
Nathan Dunfield \cite{D2} points out that this manifold fibers over
$S^1$ with orientation-reversing monodromy, with fiber a surface
of genus $2$.
\end{remark}

\begin{problem} Find the lowest volume closed Haken $3$-manifolds.\end{problem}

\begin{remark}
As of this writing the smallest known Haken manifold is the manifold obtained by $(\frac{14}{3},\frac{3}{2})$
Dehn filling on the Whitehead link complement.  It was discovered by Nathan Dunfield \cite{D1} and has volume volume approximately $2.207$.
\end{remark}

\begin{problem} Find the lowest volume closed orientable fibered
$3$-manifolds.\end{problem}

\begin{remark}  The smallest known closed fibered manifold is $m011(9,13)=v0073(3,-1)$.
It has volume $2.7317\cdots$ and is a genus-$5$ bundle. It was discovered by Saul Schleimer
and described to the authors by Nathan Dunfield \cite{D2}. \end{remark}

\begin{problem}  Find the lowest volume closed orientable $3$-manifolds with
  $\beta_1=n$.\end{problem}

\begin{remark}  Manifolds with $\beta_1=0$ are orientable, hence the Weeks manifold
is the smallest manifold, orientable or not, with $n=0$ \cite{GMM3}.  For $n=1$, the manifold
of the previous remark is the smallest known manifold.  It experimentally minimizes
all longitudinal fillings of the 1-cusped census manifolds and minimizes fillings of
longitudinal surgery on knots of $13$ crossings or less; see \cite {D2}.  For $n=2$ the
smallest known example is $(0/1,0/1)$-surgery on link $9^2_4$ also known as
$v1539(5,1)$.  This manifold, discovered by Dunfield, has volume
$4.7135\cdots$ and is a genus-$2$ fiber bundle.  Experimental work of Dunfield
\cite{D2} shows that this minimizes volume among fillings of $2$-cusped census
manifolds and $0$-surgery on homologically split $2$-component links with at most
$14$ crossings.
\end{remark}

\begin{problem} Determine the first infinite stem of closed and/or
cusped hyperbolic 3-manifolds.\label{prob:stem}\end{problem}

\begin{remark} Agol has shown in \cite{A2} that the complements of
  the Whitehead link and the pretzel link $(-2,3,8)$ are the two
lowest volume $2$-cusped hyperbolic $3$-manifolds. Their volumes are
$3.663$\ldots\ .  Thus we need to determine all the $1$-cusped
hyperbolic $3$-manifolds with volume at most $3.663$\ldots\
.\end{remark}

\begin{problem} Find the n-cusped complete finite volume hyperbolic
  3-mani\-folds of least volume. \end{problem}

\begin{remark}  As previously mention this was solved in the cases n=1,2 respectively by Cao - Meyerhoff \cite{CM} and Agol \cite{A2}.\end{remark}

\begin{problem} Develop a theory of low-volume manifolds via minimal
surfaces. \end{problem}

\begin{remark} For example if $N$ is non\-orientable, then $N$ contains an embedded
$\pi_1$-injective surface.  Such a surface can be represented by a
stable minimal surface $T$.  By Uhlenbeck, $\Area(T_0)\ge
-\pi(\chi(T_0))\ge \pi$; see \cite{Ha}.  Thus thickening and expanding the
surface creates lots of volume unless it rapidly crashes
into itself, i.e. encounters handles of index $\ge 1$.  A careful analysis
should yield lower bounds on volumes and a suitable understanding of
the geometry should yield a finite set of manifolds which contains the
lowest volume one.

If $N$ is orientable, then as announced by Pitts - Rubinstein
in \cite{PR} one can find an index-$1$ minimal Heegaard surface or a
stable $1$-sided Heegaard surface $T$.  Again careful estimates should
yield lower bound estimates on volume. Note that if $T$ is
a $1$-sided Heegaard splitting, then $N$ will also contain an
index-$1$ minimal surface disjoint from $T$.\end{remark}

\begin{Mom4Conjecture} The collection of hyperbolic Mom-$4$
manifolds is exactly the set of $117$ manifolds enumerated in
\cite{GMM2}.
\end{Mom4Conjecture}

\begin{problem}
Let $N$ be a hyperbolic $3$-manifold possessing an internal Mom-$n$
structure. Does $N$ necessarily possess a geometric Mom-$n$ structure,
i.e. a structure where the $1$-handles and $2$-handles of $\Delta$ are
neighborhoods of geodesic arcs and totally geodesic surfaces
respectively?\label{prob:interal_geometric}
\end{problem}

\begin{problem}
Let $N$ be as above. Does $N$ necessarily possess a geometric Mom-$n$
structure where the $1$-handles and $2$-handles of $\Delta$ are
neighborhoods of edges and faces in a canonical ideal cellulation
of $N$?
\end{problem}

\begin{problem} \label{internal to full} Let $N$ be a hyperbolic
$3$-manifold possessing a full internal Mom-$n$ structure.  Does $N$
possess an hyperbolic full internal Mom-$k$ structure for some $k\le
n$?\end{problem}

\begin{problem} \label{general based}Let $N$ be a hyperbolic
$3$-manifold possessing a general based internal Mom-$m$ structure.
Does $N$ posses an internal full hyperbolic Mom-$k$ structure for some
$k\le m$?  \end{problem}

\begin{remarks}
\begin{enumerate}
\item We conjecture that the answer to each of the last four problems is
``no'' for $n$ sufficiently large.

\item See \cite{GMM2} for the definition of a general based structure.
Such a structure may arise if the internal Mom structure has annular
lakes, the simplest example being the figure-$8$ knot as shown in
Figure \ref{fig:m004_badmom}.

\item A positive solution to Problem \ref{general based} for $n=2$ is
given in Lemma 4.5 of \cite{GMM2}.  A positive solution for $n=3$ under the
additional hypothesis that the Mom-$3$ structure is geometric with the
base torus cutting off a cusp is given in \cite{GMM3}.

\item A positive solution to Problem \ref{internal to full} for $n=4$ is
given in Theorem 4.1 of \cite{GMM2}.

\item Versions of the Mom-$4$ Conjecture,
Problem \ref{internal to full}, and Problem \ref{general based} are needed to
extend the Mom-$3$ technology of \cite{GMM3} to
Mom-$n$ technology.  Developing a suitable
enumeration of general based Mom-$n$ manifolds would enable one
to get around these issues.

\end{enumerate}
\end{remarks}

A striking application of the Mom technology in combination with other
geometric and topological arguments is the Lackenby-Meyerhoff
solution \cite{LM} of the long-standing Gordon conjecture.

\begin{theorem} (Lackenby - Meyerhoff)
Let $M$ be a compact 3-manifold with boundary a torus, with
interior admitting a complete finite-volume hyperbolic structure. Then
the number of non-hyperbolic Dehn fillings on $M$ is at most $10$.
\end{theorem}


\begin {problem} \label{small cusps} Find all $1$-cusped hyperbolic
$3$-manifolds with low-volume maximal cusps. In particular enumerate all
such manifolds with cusp volume at most $2.5$.\end{problem}

Since the volume of a cusp is one half the area of its boundary torus
this question is equivalent to enumerating all $1$-cusp manifolds with
maximal tori of area at most $5$.

A solution to this problem would have the following two applications.
First, much of the work in \cite{GMM3} is involved with getting
estimates on the volume of the maximal cusp. Hence a solution to
Problem \ref{small cusps} would provide a significant weapon for
attacking Problem \ref{prob:stem}.

Second, \cite{LM} makes vital use of cusp area bounds and thus Problem
\ref{small cusps} would provide an avenue towards establishing the
strong form of the Gordon Conjecture:

\begin{conjecture}The figure-$8$ knot complement is the unique
  $1$-cusped hyperbolic $3$-manifold which realizes the maximal number
  of non-hyperbolic Dehn fillings.\label{conj:Gordon}
\end{conjecture}

Of course, the dichotomy inherent in Mom technology makes it a natural
tool for working on generalizations of Conjecture
\ref{conj:Gordon}. For example:

\begin{problem}List all $1$-cusped hyperbolic $3$-manifolds that have
  $6$ or more non-hyperbolic Dehn fillings.\label{prob:prob31}
\end{problem}

There are infinitely many $1$-cusped manifolds arising from filling one
component of the Whitehead link or its sister that  have six
exceptional surgeries; see \cite{Gor}. Thus a satisfactory solution to
Problem \ref{prob:prob31} would be to list finitely many
multi-cusped manifolds
(preferably $2$-cusped or $3$-cusped) and all their fillings that yield a
$1$-cusped manifold that has six or more exceptional fillings.

\bigskip
\textbf{The following problems and discussion on number-theoretic issues were
generously provided by Walter Neumann and Alan Reid.}

\begin{problem}Suppose that $M$ is a finite-volume hyperbolic 3-manifold.
Is the volume of $M$ irrational?
\end{problem}

There seems to be no explicit reference for this ``folklore'' question.
It is worth remarking that by Apery's proof that $\zeta(3)$ is
irrational it is known that there are finite volume hyperbolic $5$-manifolds with
irrational volumes (see \cite{Ker}).

In fact, there are much more
far-reaching questions than this, namely explicit conjectures about
when volumes are linearly dependent over $\BQ$, and the same for
Chern Simons (which can be rational) \cite{NY}. This is discussed
briefly below.

As described in \cite{N}, these conjectures have appeared in different
forms in the literature. For volume they are equivalent to the
sufficiency of the Dehn invariant conjecture for $\BH^3$ scissors
congruence.

\bigskip
The following conjecture due to Milnor appears in \cite{Mil}; see also \cite{N}
and \cite{NY}. Here $\text{\cy{L}}(\theta)$ denotes the Lobachevsky function.

\begin{conjecture} (Milnor, \cite{Mil}) If we consider only angles
$\theta$ which are rational multiples of $\pi$, then every $\BQ$-linear
relation
\[
q_1 \text{\cy{L}}(\theta_1)+\cdots +q_n \text{\cy{L}}(\theta_n)=0
\]
is a consequence of the relations
\begin{eqnarray*}
\text{\cy{L}}(\pi+\theta) &=& \text{\cy{L}}(\theta),\\
\text{\cy{L}}(-\theta) &=& -\text{\cy{L}}(\theta), \\
\text{\cy{L}}(n\theta) &=& n \sum_{k\!\!\!\!\mod n} \text{\cy{L}}(\theta+k\pi/n).
\end{eqnarray*}
\label{conj:milnor}
\end{conjecture}
 

In a similar vein, the following appears as question $23$ of \cite{Th}.

\begin{problem}(Thurston, \cite{Th}) Show that the volumes of hyperbolic
$3$-manifolds are not all rationally related.
\label{prob:10_34}
\end{problem}

By work of Borel \cite{Bo}, it is known that 
given an arithmetic hyperbolic 3-manifold $M$ 
 with invariant trace-field $k$ there is a real
number $v_k$ such that $\Vol(M)$ is an integral multiple of $v_k$.

This has the practical application that it can be used to prove that an approximation to
the volume of an arithmetic manifold can be made exact. For example,
in \cite{JR} it is proved (see Lemma 3.2 and the proof of Theorem 3.1
of \cite{JR}) that $\Vol 3$ has volume $v_0$ (the volume of the
regular ideal tetrahedron in $\BH^3$).

For hyperbolic 3-manifolds, arithmetic or otherwise,
another result of Borel \cite{Bo2} (see also \cite{NY}) says:
\begin{theorem}For any number field $k$ with $r$ complex
places, there are real numbers $v_1$, $v_2$,\ldots $v_r$ such that for any
finite-volume hyperbolic $3$-manifold $M$ whose invariant trace-field is $k$
there are $r$ integers $\alpha_1$, \ldots, $\alpha_r$ such that
\[
\Vol(M) = \alpha_1v_1 + \alpha_2v_2 + \cdots + \alpha_rv_r.
\]
\label{thrm:10_35}
\end{theorem}

\begin{problem}(Neumann-Reid)
For a number field $k$ as above, identify the
real numbers $v_1$, \ldots ,$v_r$.
\end{problem}

Even when $r=1$, so that we are in the situation of the invariant trace-field
of an arithmetic hyperbolic 3-manifold, there could be non-arithmetic
hyperbolic 3-manifolds with the same invariant trace-field.

Given this, and the discussion for the arithmetic case,
some basic questions arise: For example, if one knows the
invariant trace-field of a hyperbolic manifold is $\BQ(\sqrt{-3})$,
is the volume then an integer multiple of $v_0$? This would be implied by
the Lichtenbaum conjecture, which would
more generally imply a best value for $v_1$ in Theorem \ref{thrm:10_35} whenever
$r=1$. There appears to be no good reference for this, but Gangl has
one in process \cite{Gangl} and \cite{Grayson} is relevant to this
particular case.

\bibliographystyle{amsalpha}

\end{document}